\documentclass{ifacconf}
\usepackage{amsmath}
\usepackage{amssymb}
\usepackage[x11names]{xcolor}
\usepackage{units}

\usepackage{pgfplots}
\pgfplotsset{compat=newest}
\usepackage{pgfmath,pgffor,pgfplotstable}
\usetikzlibrary{calc,decorations.markings,intersections}
\usepgfplotslibrary{fillbetween,patchplots,colorbrewer}

\usepackage{tikzscale}
\usetikzlibrary{external}
\newtheorem{assumption}{Assumption}

\newtheorem{lemma}{Lemma}
\newtheorem{theorem}{Theorem}
\newtheorem{corollary}{Corollary}
\newtheorem{definition}{Definition}
\newtheorem{remark}{Remark}

\newcommand{\nv}{{\sf v}}

\newcommand{\nf}{{\sf f}}
\newcommand{\tr}{\intercal}

\usepackage{graphicx}      % include this line if your document contains figures
\usepackage{natbib}        % required for bibliography
\begin{document}
\begin{frontmatter}

%% HOW ABOUT?
%\title{Redefining Nonlinear Feedback Design \\ Via Convex Control Lyapunov Functions\thanksref{footnoteinfo}}
%% WHAT DO YOU THINK ABOUT:
\title{Feedback Synthesis for Nonlinear Systems Via Convex Control Lyapunov Functions\thanksref{footnoteinfo}}
% \title{Convex Control Lyapunov Functions and Feedback Design for Nonlinear Systems\thanksref{footnoteinfo}}

\thanks[footnoteinfo]{RP acknowledges the support by the European Commission under the grant scheme NextGenerationEU project no. 09I01-03-V05-00002. JO acknowledges the support of the Scientific Grant Agency SR under the grant 1/0339/26, and the Slovak Research and Development Agency APVV-24-0007. Funded by the EU NextGenerationEU through the Recovery and Resilience Plan for Slovakia under the project No. 09I01-03-V04-00024.\\
Corresponding author: Boris Houska (borish@shanghaitech.edu.cn).
}
\author[First]{Mario Eduardo Villanueva}
\author[Second]{Juraj Oravec}
\author[Second]{Radoslav Paulen}
\author[Third]{Boris Houska}

\address[First]{IMT School for Advanced Studies Lucca, Italy}
\address[Second]{Slovak University of Technology in Bratislava, 
   Bratislava, Slovakia}
\address[Third]{ShanghaiTech University, China}

\begin{abstract}                % Abstract of 50--100 words
% This paper introduces computationally efficient methods for synthesizing explicit piecewise affine (PWA) feedback laws for nonlinear discrete-time systems, with guaranteed robustness and performance. The approach constructs PWA approximations to the solutions of ergodic- and infinite-horizon min-max Hamilton-Jacobi-Bellman equations by solving a tractable nonlinear program offline. The resulting feedback laws exhibit configurable storage complexity and evaluation time, determined a priori by selecting a suitable configuration template offline. The effectiveness of this framework is demonstrated on a constrained Van der Pol oscillator, where an explicit piecewise-affine controller is synthesized over a large domain with certified ergodic performance and pre-specified complexity.
%This paper introduces computationally efficient methods for synthesizing explicit piecewise affine (PWA) feedback laws for nonlinear discrete-time systems with guaranteed robustness and performance. The approach constructs PWA approximations to the solutions of infinite-horizon min–max Hamilton–Jacobi–Bellman equations by solving a tractable nonlinear program offline, yielding feedback laws with configurable storage complexity and evaluation time determined a priori via a chosen configuration template. The effectiveness of the framework is illustrated on a constrained Van der Pol oscillator, where an explicit PWA controller with certified ergodic performance and pre-specified complexity is synthesized over a large domain.
This paper introduces computationally efficient methods for synthesizing explicit piecewise affine (PWA) feedback laws for nonlinear discrete-time systems, ensuring robustness and performance guarantees. The approach proceeds by optimizing a configuration-constrained PWA approximation of the value function of an infinite-horizon min–max Hamilton–Jacobi–Bellman equation. Here, robustness and performance are maintained by enforcing the PWA approximation to be a generalized control Lyapunov function for the given nonlinear system. This enables the generation of feedback laws with configurable storage complexity and pre-determined evaluation times, based on a selected configuration template. The framework's effectiveness is demonstrated through a constrained Van der Pol oscillator case study, where an explicit PWA controller with certified ergodic performance and specified complexity is synthesized over a large operational domain.
\end{abstract}

\begin{keyword}
Non-Linear Control Systems,
Optimal Control,
Robust Control,
Control Lyapunov Functions,
Explicit Model Predictive Control
\end{keyword}

\end{frontmatter}
%===============================================================================

\section{Introduction}

% Control Lyapunov Functions (CLFs) were originally introduced by \cite{Zubov1965} and further developed by \cite{Artstein1983} and \cite{Primbs1999}. Over time, CLFs have emerged as a cornerstone in the stability analysis and controller synthesis for both linear and nonlinear systems. This prominence is attributed to their ability to provide a systematic framework for designing feedback controllers that ensure asymptotic stability, even in the presence of constraints and uncertainties. A comprehensive survey by \cite{Giesl2015} highlights the evolution of CLF methodologies, emphasizing their maturity in addressing modern control challenges.

Control Lyapunov Functions (CLFs), introduced by \cite{Zubov1965} and developed in \cite{Artstein1983,Primbs1999}, have become a cornerstone of stability analysis and control synthesis for both linear and nonlinear control systems. They provide a systematic framework for designing feedback controllers that ensure asymptotic stability, even in the presence of constraints and uncertainties. A comprehensive survey by \cite{Giesl2015} traces the evolution of CLF methodologies and highlights their maturity in addressing modern control challenges.

% General methods for computing or approximating the maximal control invariant domain of CLF functions have been developed within the set-valued analysis literature, and we refer to \cite{Aubin2009} for an overview. In particular, polyhedral computing methods have gained prominence due to their efficiency in handling constrained systems, as discussed in \cite{Blanchini2008} and \cite{Fukuda2020}; see also~\cite{Houska2025} for a recent survey. These polyhedral computing methods have also been extended for constructing piecewise affine Lyapunov functions. For instance, Blanchini's pioneering work \citep{Blanchini1999a} demonstrated how special classes of polyhedral gauge functions can serve as CLF candidates, while subsequent studies have explored non-quadratic CLFs \citep{Blanchini1995}. Recent extensions of this work can be found in \cite{Rakovic2020} and \cite{Rakovic2024}, where special classes of Minkowski-Lyapunov functions are used for constructing stablilizing but sub-optimal controllers. Moreover, \cite{Houska2025a} presents a convex optimization-based method for computing configuration-constrained polyhedral CLFs, although these methods are currently limited to constrained linear systems.

Methods for computing or approximating the maximal control invariant domain of CLF functions have been developed within the set-valued analysis community and we refer to~\cite{Aubin2009} for an overview. Among these, polyhedral computing has become prominent for its efficiency in handling constraints, as discussed in \cite{Blanchini2008}, \cite{Fukuda2020}, and recently in~\cite{Houska2025}. Beyond invariant set computation, polyhedral techniques have also been used to construct piecewise affine CLFs. \cite{Blanchini1999a} pioneered the use of polyhedral gauge functions as CLF candidates and related studies have explored non-quadratic CLFs~\citep{Blanchini1995}. \cite{Rakovic2020} and~\cite{Rakovic2024} extended this work to special classes of Minkowski-Lyapunov functions for designing stabilizing, generally sub-optimal, controllers. Additionally, in recent years, a variety of methods for computing piecewise affine CLFs by using neural networks have been proposed. We refer to~\cite{Fabiani2022}, \cite{Gruene2025}, and \cite{He2024} for the latest developments and an overview of these methods. Moreover, a convex optimization-based approach for computing configuration-constrained polyhedral CLFs for constrained linear systems was introduced in~\cite{Houska2025a}.

% CLFs are commonly featured as terminal cost functions in Model Predictive Control (MPC) for both linear and nonlinear systems, ensuring recursive feasibility and stability \citep{Mayne2000, Rawlings2015}. This principle extends to MPC formulations for uncertain systems, particularly in min-max robust MPC frameworks, as demonstrated in studies by \cite{Bemporad2003}, \cite{Diehl2004}, and \cite{Kerrigan2004}, all of which build upon the foundational work of \cite{Witsenhausen1968} on min-max dynamic programming recursions. Additionally, alternative approaches for constructing CLFs in uncertain nonlinear systems can be found in the Tube MPC literature. For example, \cite{Villanueva2017} proposes a method where nonlinear terms are over-approximated by treating them as uncertainties.

CLFs are widely used as terminal cost functions in Model Predictive Control (MPC) for linear and nonlinear systems to guarantee recursive feasibility and stability \citep{Mayne2000,Rawlings2015}. This idea extends to MPC under uncertainty, particularly in min–max MPC frameworks \citep{Bemporad2003,Diehl2004,Kerrigan2004}, which build on Witsenhausen's seminal work on min–max dynamic programming~\citep{Witsenhausen1968}. Further approaches to CLF construction for uncertain nonlinear systems arise in the Tube MPC literature, where nonlinear terms are treated as uncertainties~\citep{Villanueva2017}.

% Beyond their applications in MPC, CLFs can also be used to construct explicit feedback laws. This property is leveraged in exact explicit MPC methods, as introduced by~\cite{Bemporad2002} and implemented in software tools such as MPT~\citep{Herceg2013}. While these exact methods are typically restricted to linear systems and their associated piecewise quadratic CLFs, recent advancements have been made. For instance, \cite{Pappas2021,Pappas2021_explicit} presents methods for computing explicit solutions to convex and non-convex quadratically constrained quadratic programming (QCQP) problems. Additionally, polytopic approximation techniques in explicit MPC are discussed in \cite{Jones2010}. \textcolor{blue}{TODO: review more related article by Stratos, Tor Arne Johansen, Bart De Schutter, etc.}, \cite{Johansen2004}, \cite{Ganguly2025}

Beyond their applications in MPC, CLFs also serve as a basis for constructing explicit feedback laws. This idea is exploited in exact explicit MPC methods, as introduced by~\cite{Bemporad2002} and implemented in software tools such as MPT~\citep{Herceg2013}. While these exact formulations are typically limited to linear systems with piecewise quadratic CLFs, several recent extensions have been proposed. For instance, \cite{Pappas2021} developed methods for computing explicit solutions to convex and non-convex quadratically constrained quadratic programs (QCQPs) and polytopic approximations for explicit MPC were discussed in~\cite{Jones2010}. Orthogonally, \cite{Ganguly2025} proposed an interpolation approach for approximating CLFs on a grid.

% Finally, under certain assumptions, the value function of infinite-horizon optimal control problems can also serve as a CLF for the underlying system (see \cite{Rawlings2015, Giesl2015}). Consequently, any dynamic programming method can, in principle, be used to construct a CLF \citep{Bellman1957}. However, it should be noted that solving the associated Hamilton-Jacobi-Bellman (HJB) partial differential equation is computationally demanding, particularly for high-dimensional systems. As a result, many numerical approaches, such as those developed in \cite{Luss2000} and \cite{Gruene2002}, are typically restricted to low-dimensional state spaces. In this context, it is also important to acknowledge recent operator-theoretic frameworks for nonlinear control, including those based on Koopman operators \citep{Korda2018} or more general Fokker-Planck-Kolmogorov operators, as used in modern operator-theoretic dynamic programming \citep{Houska2025c}. Most of these dynamic programming methods focus on computing accurate approximations of the value function for infinite-horizon optimal control problems to derive explicit, approximately optimal control laws. However, these approximations of the optimal value function are not necessarily CLFs, and thus, the corresponding explicit feedback laws may not provide strong stability guarantees.

Under suitable assumptions, the value function of an infinite-horizon optimal control problem (OCP) can itself serve as a CLF for the underlying system \citep{Rawlings2015,Giesl2015}, so that any dynamic programming method can, in principle, be used to construct a CLF~\citep{Bellman1957}. In practice, however, solving the associated Hamilton–Jacobi–Bellman (HJB) partial differential equation is computationally demanding, especially in high dimensions~\citep{Bardi1997}. Therefore, many numerical approaches are restricted to low-dimensional state spaces~\citep{Luss2000,Gruene2002}. Related developments include operator-theoretic frameworks for nonlinear control, such as those based on Koopman operators~\citep{Korda2018} or more general Fokker–Planck–Kolmogorov operators in modern operator-theoretic dynamic programming~\citep{Houska2025c}. Most of these methods aim to approximate the infinite-horizon value function to obtain explicit, approximately optimal feedback laws. However, the resulting value function approximations are not necessarily CLFs, and the corresponding feedback laws therefore do not always provide strong stability guarantees.

\subsection*{Contributions}

This paper introduces a control synthesis framework for constrained, uncertain nonlinear systems with economic stage costs, leveraging the theoretical rigor of CLFs. All computationally demanding steps are carried out offline, eliminating the need for online optimization while retaining predefined closed-loop performance. In contrast to traditional Explicit MPC, however, the approach yields nonlinear control laws with a predefined structural complexity, leading to predictable memory requirements and implementation costs. The explicit control policies admit deterministic worst-case runtimes, ensuring real-time feasibility. Their explicit form facilitates rigorous a priori certification of closed-loop properties, including robustness and performance.

\subsection*{Outline}

% \begin{itemize}
%     \item The main theoretical results are summarized in Lemma~\ref{lem::ergodic} and Theorem~\ref{thm::CLF}.
%     \item The optimization problem for designing generalized CLFs is summarized in~\eqref{eq::DESIGN_PROBLEM}.
% \end{itemize}
%The remainder of the paper is organized as follows.

Section~\ref{sec::OC} introduces the general nonlinear discrete-time system model together with the infinite-horizon optimal control formulation and the definition of generalized control Lyapunov functions. Section~\ref{sec::CCCLF} presents the main theoretical developments and results, summarized in Theorem~\ref{thm::CLF}, and formulates the convex optimization problem for CLF design in \eqref{eq::DESIGN_PROBLEM}. Section~\ref{sec::case_studies} demonstrates the applicability of the proposed framework on a nonlinear Van der Pol oscillator. Section~\ref{sec::conclusions} concludes the paper.

\subsection*{Notation}
We denote the convex hull by $\mathrm{cvx}$.
The unit simplex $\Delta_n \subset \mathbb{R}^{n}$ with $n = n_{\mathrm{x}}+1$ is denoted by
\[
\Delta_n \ = \ \left\{ \ \theta \in \mathbb R^{n} \ \middle| \ \theta \geq 0, \ \sum_{i=1}^n \theta_i = 1 \ \right\}.
\]
The indicator function of a set $\mathbb X \subseteq \mathbb R^{n_{\mathrm{x}}}$ 
is given by
\[
I_{\mathbb X}(x) \ = \ \left\{
\begin{array}{ll}
0 & \text{if} \ x \in \mathbb X \\
\infty & \text{otherwise}.
\end{array}
\right.
\]
More broadly, throughout this paper, we extensively use extended value conventions, wherein the objective value of a constrained minimization problem is defined as $\infty$ whenever the problem is infeasible.

\section{Optimal Control}
\label{sec::OC}
This section focuses on designing feedback laws for time-autonomous discrete-time systems. We begin by introducing the system model and its properties, followed by a discussion on infinite-horizon optimal control and the role of CLFs.

\subsection{Nonlinear System Model}
We consider nonlinear discrete-time systems of the form
\begin{equation}
\label{eq::sys}
\forall k \in \mathbb N, \qquad x_{k+1} = f(x_k,u_k) + w_k.
\end{equation}
Here, $x_k \in \mathbb R^{n_{\mathrm{x}}}$, $u_k \in \mathbb R^{n_{\mathrm{u}}}$, and $w_k \in \mathbb R^{n_{\mathrm{x}}}$ denote, respectively, state, control, and disturbance vectors at time $k \in \mathbb N$. Similarly, $\mathbb X \subseteq \mathbb R^{n_{\mathrm{x}}}$, $\mathbb U \subseteq \mathbb R^{n_{\mathrm{u}}}$, and $\mathbb W \subseteq \mathbb R^{n_{\mathrm{x}}}$ denote given convex state, control, and disturbance sets. Additionally, $\mathbb X$ is assumed closed, and $\mathbb U$ and $\mathbb W$ compact. The function $f: \mathbb R^{n_{\mathrm{x}}} \times \mathbb R^{n_{\mathrm{u}}} \to \mathbb R^{n_{\mathrm{x}}}$ is required to satisfy the following Lipschitz continuity assumption.

\begin{assumption}
\label{ass::f}
There exist constants $\gamma,\alpha \geq 0$ such that
\begin{align*}
& \max_{\theta \in \Delta_n} \ \left\| \sum_{i=1}^{n} \theta_i f(v_i,u_i) - f\left( \sum_{i=1}^{n} \theta_i v_i, \sum_{i=1}^{n} \theta_i u_i \right) \right\|_\infty \\
& \leq \ \gamma \cdot \max_{i,j} \left\|
\left[
\begin{array}{c}
v_i - v_j \\
u_i - u_j
\end{array}
\right]
\right\|^\alpha.
\end{align*}
holds for all $v_1,\ldots,v_n \in \mathbb X$ and all $u_1,\ldots,u_n \in \mathbb U$. 
\end{assumption}
If $f$ is Lipschitz continuous, Assumption~\ref{ass::f} is trivially satisfied for $\alpha = 1$. For functions with bounded second-order derivatives, we set $\alpha = 2$. If $f$ is linear, the assumption is satisfied with $\gamma = 0$, motivating the interpretation of $\gamma$ as a bound on the nonlinearity of $f$.

\subsection{Infinite-Horizon Optimal Control}

In practice, a non-negative and Lipschitz continuous stage cost function $L: \mathbb{X} \times \mathbb{U} \to \mathbb{R}_+$ is often given.\footnote{If $L$ is merely bounded from below, a constant offset can be added to ensure non-negativity without altering the solution of its associated optimal control problem.} The sequence of optimal value functions is defined by the reverse-indexed dynamic programming recursion
% $$I_{\mathbb X}(x) = \begin{cases} 0 &\text{if} \ x\in\mathbb{X}\\
% 0 &\text{otherwise}\end{cases}$$
\begin{equation}\label{eq::HJB} 
\!\!\!\!
\begin{aligned}
J_{k+1}(x) &= \min_{u \in \mathbb U}  \max_{w \in  \mathbb W}  L(x,u) + I_{\mathbb X}(x) + J_{k}(f(x, u)+w) \\
J_{0}(x) &= I_{\mathbb X}(x).
\end{aligned}
\end{equation}
This recursion remains well-defined and returns non-negative and lower semi-continuous functions $J_k \geq 0$ for all integers $k \in \mathbb N$.
%Here, we set $J_k(x) = \infty$ whenever its associated OCP is infeasible.
This statement holds because $f$ and $L$ are Lipschitz continuous, while $\mathbb X$, $\mathbb U$, and $\mathbb W$ are all convex and closed. In contrast, the existence of a limit
\[
J_\infty(x,d) \ = \ \lim_{k \to \infty} \, [ J_k(x) - k \cdot d ],
\]
for a drift $d \in \mathbb R$, is much harder to establish. If it exists, $(J_\infty,d)$ satisfies the infinite horizon Bellman equation
\begin{equation}
\begin{aligned}
& J_\infty(x,d) + d \\
\label{eq::infHJB}
&\quad=  \min_{u \in \mathbb U}  \max_{w \in  \mathbb W} L(x,u) + I_{\mathbb X}(x) + J_\infty(f(x,u)+w,d).
\end{aligned}
\end{equation}
This follows, at least formally, by taking the limit for $k \to \infty$ in~\eqref{eq::HJB}, see~\cite{Houska2025c} and~\cite{Namah1999}. Nevertheless, it should be kept in mind that---at least in the general setting for arbitrary nonlinear functions $f$ and $L$---nothing is known about the existence and uniqueness of~\eqref{eq::infHJB}.

\begin{assumption}
\label{ass::L}
% The function $L$ is Lipschitz continuous and non-negative. Moreover, t
There exist constants $\sigma \geq 0$ and $\beta \geq 1$ such that
\begin{align*}
& \max_{\theta \in \Delta_n} \ \left( L\left( \sum_{i=1}^{n} \theta_i v_i, \sum_{i=1}^{n} \theta_i u_i \right) - \sum_{i=1}^{n} \theta_i L(v_i,u_i) \right) \\
& \leq \ \sigma \cdot \max_{i,j} \left\|
\left[
\begin{array}{c}
v_i - v_j \\
u_i - u_j
\end{array}
\right]
\right\|^\beta.
\end{align*}
holds for all $v_1,\ldots,v_n \in \mathbb X$ and all $u_1,\ldots,u_n \in \mathbb U$.
\end{assumption}

This assumption is satisfied with $\beta = 1$ for all non-negative and Lipschitz continuous $L$. If $L$ is twice differentiable and the eigenvalues of its Hessian are bounded from below, we may set $\beta = 2$. If $L$ is convex, we may set $\sigma = 0$, motivating the interpretation of $\sigma$ as a measure of non-convexity.

\subsection{Control Lyapunov Functions}
\label{sec::clf}

Solving Bellman's equation \eqref{eq::infHJB} for infinite-horizon optimal control is often computationally intractable. To address this challenge, we introduce a relaxation framework based on generalized control Lyapunov functions (CLFs).

\begin{definition}%[Generalized Control Lyapunov Function]
\label{def::CLF}
For a given drift parameter $d \in \mathbb{R}$, a convex function $M: \mathbb{R}^{n_{\mathrm{x}}} \to \mathbb{R}_+ \cup \{\infty\}$ is called a generalized CLF for the supply rate $(L - d)$ if it is non-negative and radially unbounded and satisfies
the dissipation inequality 
\begin{equation}\label{eq::CLF}
M(x) + d \geq \min_{u \in \mathbb{U}} \max_{w \in \mathbb{W}} L(x,u) + I_{\mathbb{X}}(x) + M(f(x,u) + w) 
\end{equation}
for all $x \in \mathbb{R}^{n_{\mathrm{x}}}$.
\end{definition}
The continuity of $f$ (see Assumption~\ref{ass::f}) and lower semi-continuity of $L$ (see Assumption~\ref{ass::L}), together with compactness of $\mathbb{U}$ and $\mathbb{W}$, ensure that the min-max operations in~\eqref{eq::CLF} are well-defined in the extended value sense. Since~\eqref{eq::CLF} is invariant under constant offsets to $M$, we may assume without loss of generality that
\[
\min_{x \in \mathbb{R}^{n_{\mathrm{x}}}} M(x) = 0\;.
\]
Moreover, the minimizer exists, because $M$ is convex and radially unbounded.

If the supply rate $(L - d)$ is non-negative on $\mathbb{X} \times \mathbb{U}$, Definition \ref{def::CLF} reduces to the standard CLF condition found in the literature; see \cite{Giesl2015}, \cite{Rawlings2015}, and~\cite{Houska2025c}. However, the generalized formulation remains useful even when $(L - d)$ is not positive definite, as the ergodic performance bound established in the following statement holds.

\begin{lemma}%[Ergodic Performance Bound]
\label{lem::ergodic}
Let Assumptions \ref{ass::f} and \ref{ass::L} hold, $d \in \mathbb R$ be a given drift parameter with $M$ being an associated generalized CLF, and let $(J_k)_{k \in \mathbb{N}}$ be the cost sequence defined by \eqref{eq::HJB}. Then, for any $x \in \mathbb{X}$ with $M(x) < \infty$, the ergodic limit
\begin{align}
\label{eq::dinfty}
d_{\infty}(x) = \limsup_{k \to \infty} \frac{J_k(x)}{k+1}
\end{align}
is bounded and satisfies $d_{\infty}(x) \leq d$.
\end{lemma}

\begin{pf}
Our first goal is to show that
\begin{align}
\label{eq::ineq55}
M(x) \geq J_k(x) - (k+1) \cdot d.
\end{align}
We proceed by induction. For the induction start, we note 
$L$ and $M$ are non-negative and $d < \infty$. Thus, 
\eqref{eq::CLF} implies
\[
M(x)+d \geq I_{\mathbb{X}}(x) = J_0(x) \quad \Longrightarrow \quad M(x) \geq J_0(x)-d.
\]
Next, assume that~\eqref{eq::ineq55} holds for a given $k \in \mathbb{N}$. Using the dynamic programming principle \eqref{eq::HJB} and the CLF condition \eqref{eq::CLF}, we obtain
\[
M(x) \geq J_{k+1}(x) - (k+2) \cdot d,
\]
confirming the induction step. Consequently, we have
\begin{align}
\label{eq::dinfty_bound}
\frac{J_k(x)}{k+1} \leq d + \frac{M(x)}{k+1}.
\end{align}
For $x$ with $M(x) < \infty$, the sequence $J_k(x)/(k+1)$ is bounded. By taking the limit superior as $k \to \infty$ in \eqref{eq::dinfty_bound} yields $d_{\infty}(x) \leq d$.
 \hfill$\diamond$
 \end{pf}

\begin{remark}
As discussed above, for $d > 0$, $M$ does not constitute a CLF in the classical sense (which requires the supply rate to be non-negative). However, as discussed in the above lemma, it can be interpreted as a generalized CLF for economic cost functions, where performance is characterized by average cost rather than stability. As we shall see in the following section, this relaxation has the advantage that one can regard $d$ as a free optimization variable, which can, for example, be minimized subject to the dissipation condition in~\eqref{eq::CLF}.
\end{remark}

\section{Configuration-Constrained Control Lyapunov Functions}
~\label{sec::CCCLF}
The generalized CLF condition~\eqref{eq::CLF} does not fully resolve the computational challenges of solving Bellman's equation~\eqref{eq::infHJB}, as the search for general CLFs remains computationally demanding. To address this problem, we propose restricting our search to functions whose epigraphs are \emph{configuration-constrained polyhedra} with pre-specified complexity. This approach leverages existing configuration-constrained polytopic computing techniques introduced in~\cite{Villanueva2024} (see also~\cite{Houska2025}). In contrast to standard piecewise affine CLFs--where complexity can grow undesirably and vertex enumeration can become computationally prohibitive--the configuration-constrained approach fixes the representation complexity in advance and parameterizes the vertices efficiently, as elaborated in the sections below.

\subsection{Configuration-Constrained Polyhedral Epigraphs}
We introduce a specific class of configuration-constrained polyhedra $\mathcal P(z) \subseteq \mathbb R^n$, with $n = n_{\mathrm{x}} + 1$, that admit a parameterized Minkowski-Weyl representation suitable for discretizing~\eqref{eq::CLF}. In this context, $\mathcal{P}(z)$ is defined as
\begin{eqnarray}
\mathcal P(z) &=&
\left\{ \
\left[
\begin{array}{c}
x \\
y
\end{array}
\right] \ \middle| \ x \in \mathbb R^{n_{\mathrm{x}}}, \ y \in \mathbb R, \
\left[
\begin{array}{cc}
G_1 & 0 \\
G_2 & h_2
\end{array}
\right] 
\left[
\begin{array}{c}
x \\
y
\end{array}
\right] \leq z \
\right\} \notag \\
\label{eq::ccond}
&=& \mathrm{cvx}
\left( 
\left[
\begin{array}{c}
V_1 z \\
s_1^\tr z
\end{array}
\right] , \ldots, 
\left[
\begin{array}{c}
V_\nv z \\
s_\nv^\tr z
\end{array}
\right]
\right) + \mathbb R_+ e_n,
\end{eqnarray}
for all parameters $z$ with $E z \leq 0$. Definition~\eqref{eq::ccond} relies on the following notation, introduced in~\citet[Sect.~3]{Villanueva2024} (see also~\citet[Sect.~3]{Houska2025a}):
\begin{enumerate}
\item The rows of $G_1 \in \mathbb R^{\nf_1 \times n_{\mathrm{x}}}$ and $[G_2 \ h_2]\in \mathbb R^{\nf_2 \times (n_{\mathrm{x}} +1)}$ define the facet normals for $\mathcal P(z)$ and are such 
that $G_1 x \leq 0$ implies $x=0$ and $h_2 < 0$;

\item the matrices $V_1,\ldots,V_{\nv} \in \mathbb R^{n_{\mathrm{x}} \times \nf}$ and vectors $s_1,\ldots,s_{\nv} \in \mathbb R^{\nf}$, with $\nf = \nf_1 + \nf_2$, are used to locate each of the $\nv$ vertices of $\mathcal P(z)$; and

\item the edge matrix $E \in \mathbb{R}^{\e \times \nf}$ encodes the combinatorial face configuration, such that~\eqref{eq::ccond} holds for a specific parameter $z \in \mathbb R^f$ if and only if $E z \leq 0$. For its construction, see e.g.~\citet[Sect. 3.5]{Villanueva2024} and~\citet[Prop.~1]{Houska2024a}.
\end{enumerate}

Now, $\mathcal P(z)$ is the epigraph of the piecewise affine function
\begin{equation}
M_z(x) = 
\begin{cases}
\underset{i \in \{1,\ldots,{\sf f}_2\}}{\max} \ \dfrac{(z_2 - G_2 x)_{i}}{(h_2)_i}  \ &\text{if} \ \ G_1 x \leq z_1 \,\\[0.1cm]
\quad \infty & \text{otherwise}.
\end{cases}
\label{eq::Mz}
\end{equation}
Below, we shall use this parameterization to discretize~\eqref{eq::CLF}.

\subsection{Configuration-Constrained Generalized CLFs}
Let $M_z$ be defined as in \eqref{eq::Mz}, with $v_i = V_i z$ for $i \in \{ 1, \ldots, \nv \}$ denoting the projection of its vertices onto the state space. 
For the given configuration and each $v_i$, let $\mathcal N_i$ be the set of indices $j$, including $i$, such that $v_j$ shares a facet with $v_i$. Moreover,
% For each vertex $v_i$, define $\mathcal{N}_i$ as the set of indices $j$ corresponding to vertices sharing a common facet with $v_i$ under the given configuration. To quantify the nonlinearity and non-convexity effects, we introduce the shorthands
\[
\mathcal E_i(u,v) \! \ = \ \! \left\{ \ \! e \in \mathbb R^{n_\mathrm{x}} \ \middle| \ \| e \|_\infty \leq \gamma \cdot  \max_{j,k \in \mathcal N_i} \left\|
\left[
\begin{array}{c}
v_j - v_k \\
u_j - u_k
\end{array}
\right]
\right\|^\alpha \
\! \right\}
\]
quantifies the nonlinearity in $f$ per Assumption~\ref{ass::f} while
\[
\kappa_i(v,u) \ = \ \sigma \cdot \max_{j,k \in \mathcal N_i} \left\|
\left[
\begin{array}{c}
v_j - v_k \\
u_j - u_k
\end{array}
\right]
\right\|^\beta,
\]
quantifies the nonconvexity in $L$ per Assumption~\ref{ass::L}. We also introduce the 
shorthands\footnote{Addition of sets is given by $\mathbb{A}+\mathbb{B} = \{ a+b \,|\, a\in \mathbb{A}, \, b\in \mathbb{B} \}$.} $\mathbb D_i = \mathcal E_i(u,v) + \mathbb W$, and $(v,u)$ to collect 
the vertices and corresponding control inputs. The following theorem establishes sufficient conditions for $M_z$ to satisfy the generalized CLF condition~\eqref{eq::CLF}\footnote{We omit the $(u,v)$ dependencies as no 
confusion should arise.}.

\begin{theorem}
\label{thm::CLF}
Let Assumptions~\ref{ass::f} and~\ref{ass::L} hold. Let $z$ satisfying $Ez \leq 0$ and $d \in \mathbb R$ be given. If the inequality
% \begin{eqnarray}
% M_z(v_i) + d \ \geq \ &\max_{w \in \mathbb W} \max_{\epsilon_i \in \mathcal  E_i(v,u)} & L(v_i,u_i) + \kappa_i(v,u) \notag \\
% \label{eq::maxMz}
% & & + M_z(f(v_i,u_i) + \epsilon_i + w) \qquad
% \end{eqnarray}
\begin{equation*}
M_{z}(v_i) + d \geq \max_{\delta \in {\mathbb D}_i} L(v_i,u_i)
+ \kappa_i + I_{\mathbb X}(v_i) + M_z(f(v_i,u_i)+\delta),
\end{equation*}
holds for all vertices $v_1,\ldots,v_{\sf v}$ and given feasible control inputs $u_1,\ldots,u_\nv \in \mathbb U$, then $(M_z,d)$ satisfies~\eqref{eq::CLF}.
\end{theorem}

\begin{pf}
Let $\overline x \in \mathbb R^{n_{\mathrm{x}}}$ be arbitrary and satisfy $M_z(\overline x) < \infty$. Due to Caratheodory's theorem, there exists an index set $\mathcal N$, with $|\mathcal N| \leq n_{\mathrm{x}}+1$, such that we have
\[
\overline x \ = \ \sum_{i \in \mathcal N} \theta_i v_i \quad \text{with} \quad \theta_i \geq 0, \ \ \sum_{i \in \mathcal N} \theta_i = 1\;.
\]
% We may assume that $\mathcal N \subseteq \mathcal N_k$ for some $k \in \{ 1, \ldots, \nv \}$, as we can choose $\mathcal N$ such that its associated vertices are all on the same facet. Define the shorthands
Moreover, $\mathcal N$ can be chosen such that $\{ v_i \mid i\in\mathcal N$ \} is a subset of one region of the polyhedral partition. Let
\[
\overline u  =  \sum_{i \in \mathcal N} \theta_i u_i, \ \ \overline \kappa = \min_{i \in \mathcal N} \kappa_i, \ \ \text{and} \ \ \overline{\mathbb D} = \bigcap_{i \in \mathcal N} \mathbb D_i.
\]
Assumption~\ref{ass::L} guarantees that the construction is such that
\begin{equation}
\label{eq::ineq1}
\sum_{i \in \mathcal N} \theta_i L(v_i,u_i) + \overline \kappa \ \geq \ L(\overline x,\overline u).
\end{equation}
Likewise, since $M_z$ is convex, Assumption~\ref{ass::f} guarantees
\begin{align}
& \sum_{i \in \mathcal N}  \max_{\delta \in \overline{\mathbb D}}\,  \theta_i \cdot M_z( f(v_i,u_i) + \delta ) \notag \\
&\overset{\text{Jensen's ineq.}}{\geq}  \max_{\delta \in \overline{\mathbb D}} \, M_z\left( \sum_{i \in \mathcal N}  \theta_i f(v_i,u_i) + \delta \right) \notag \\ 
&\overset{\text{Assumption~\ref{ass::f}}}{\geq} \max_{w \in \mathbb W}\,  M_z\left( f( \overline x, \overline u) + w \right). \label{eq::ineq2}
\end{align}
Finally, the hypothesis and the above relations yield
\begingroup
\allowdisplaybreaks
\begin{align}
&M_z(\overline{x})+d = M_z\left(\sum_{i\in\mathcal N}\theta_i v_i\right)+d = 
\sum_{i\in\mathcal N}\theta_i (M_z(v_i)+d) \notag\\
&\geq \sum_{i\in\mathcal N}\theta_i\max_{\delta \in \mathbb D_i}  L(v_i,u_i) + \kappa_i +I_{\mathbb X}(v_i) + M_z(f(v_i,u_i)+\delta) \notag\\ 
&\geq \sum_{i\in\mathcal N}\theta_i\max_{\delta \in \overline{\mathbb D}}  L(v_i,u_i) + \overline \kappa +I_{\mathbb X}(v_i) + M_z(f(v_i,u_i)+\delta) \notag\\ 
&= \sum_{i\in\mathcal N}\theta_i L(v_i,u_i) + \overline \kappa +\sum_{i\in\mathcal N}\theta_i I_{\mathbb X}(v_i) \notag\\
&\hspace{2.68cm}+ \sum_{i\in\mathcal N}\theta_i\max_{\delta \in \overline{\mathbb D}}  M_z(f(v_i,u_i)+\delta) \notag\\ 
&\geq \max_{w\in\mathbb  W} L(\overline x, \overline u) + I_\mathbb{X}(\overline x) + 
M_{z}(f(\overline x, \overline u) + w) \notag\\
&\geq \min_{u\in\mathbb U} \max_{w\in\mathbb  W} L(\overline x, u) + I_\mathbb{X}(\overline x) + 
M_{z}(f(\overline x, u) + w).\label{eq::ineq3}
\end{align}
\endgroup
% \begin{eqnarray}
% M_z(\overline x) &=& M_z\left( \sum_{i \in \mathcal N} \theta_i v_i \right) \ = \ \sum_{i \in \mathcal N} \theta_i M_z\left( v_i \right) \notag \\
% &\geq& \sum_{i \in \mathcal N} \theta_i \left( \max_{w_i \in \mathbb W} \max_{\epsilon_i \in \mathcal E_i(v,u) } \left( L\left( v_i,u_i \right) + \kappa_i(v,u) \right. \right. \notag \\
% & & \left. \left. \phantom{\max_{\epsilon_i \in \mathcal E_i(v,u) }}
% + I_{\mathbb X}(v_i) + M_z(f(v_i,u_i) + \epsilon_i + w_i ) \right) \right) \notag \\
% &\geq& \sum_{i \in \mathcal N} \theta_i \left( \max_{w_i \in \mathbb W} \max_{\epsilon_i \in \overline{\mathcal E}} \left( L\left( v_i,u_i \right) + \overline{\kappa} + I_{\mathbb X}(v_i) \right. \right. \notag \\
% & & \left. \phantom{\max_{\epsilon_i \in \overline{\mathcal E}} \max_{\epsilon_i \in \overline{\mathcal E}}} \left. + M_z(f(v_i,u_i) + \epsilon_i + w_i ) \right) \right) \notag \\
% &=& \sum_i \theta_i L(v_i,u_i) + \overline{\kappa} + \sum_i \theta_i I_{\mathbb X}(v_i) \notag \\
% & & + \sum_i \max_{w_i \in \mathbb W} \max_{\epsilon_i \in \overline{\mathcal E}} \theta_i M_z( f(v_i,u_i) + \epsilon_i + w_i ) \notag \\
% &\geq& \max_{w \in \mathbb W} \left( L\left( \overline x, \overline u \right) + I_{\mathbb X}(\overline x) + M_z( f(\overline x, \overline u) + w )  \right) \notag \\
% &\geq& \min_{\mathrm{u}} \max_{w \in \mathbb W} \left( L\left( \overline x, u \right) + I_{\mathbb X}(\overline x) + M_z( f(\overline x, u) + w )  \right). \quad \label{eq::ineq3}
% \end{eqnarray}
Here, the equations in the first line hold, after recalling that all $v_i$ with $i\in\mathcal N$ are in the same region of the domain of $M_z$---where $M_z$ is 
affine. The second line follows after substituting the vertex inequality in the 
statement of the theorem. The third line follows from $\overline \kappa \leq \kappa_i$ and 
$\overline{\mathbb D}\subseteq \mathbb D_i$. The second to last inequality follows 
from convexity of the indicator (Jensen's inequality) as well as substituting~\eqref{eq::ineq1} and~\eqref{eq::ineq2}. Since~\eqref{eq::ineq3} holds for all $\overline x$ with $M_z(\overline x)$, we conclude that $(M_z,d)$ satisfies~\eqref{eq::CLF}, as claimed. \hfill$\diamond$
\end{pf}

\subsection{Optimization of CLFs}
\label{sec::OptimalCLFs}
The conditions from Theorem~\ref{thm::CLF} can be used to design a~polyhedral CLF. Upon expansion, one can write a~nonlinear program of the form
\begin{equation}\label{eq::DESIGN_PROBLEM}
\begin{alignedat}{2}
&\min_{d,u,y,z,\lambda,\kappa} \ && \rho(z,d) \\
&\hphantom{_{d,u,}}\text{s.t.} &&  \left\{
\begin{aligned}
&\forall i \in \{ 1, \ldots, \nv \}, \forall j \in \mathcal N_i, \\
&s_i^\intercal z \geq L(V_i z,u_i) + \kappa_i - d + y_i \\
&Gf( V_i z, u_i) + hy_i + \overline w+ \lambda_i  g  \leq  z \\
&\gamma \left\|
\begin{bmatrix}
(V_i - V_j)z \\
u_i-u_j
\end{bmatrix}
\right\|^\alpha
 \leq  \lambda_i \\
&\sigma \left\|
\begin{bmatrix}
(V_i - V_j)z \\
u_i-u_j
\end{bmatrix}
\right\|^\beta
 \leq  \kappa_i \\
&E z  \leq  0, \ V_i z \in \mathbb X, \ u_i \in \mathbb U\;.
\end{aligned}
\right.
\end{alignedat}
\end{equation}
Here, $G^\intercal = [G_1^\intercal, \, G_2^\intercal]$ and $h^\intercal = [0^\intercal, \; h_2^\intercal]$, as well as 
\[
\forall j \in \{1, \ldots, \nf \} , \ \  g_j = \max_{\| n \|_\infty \leq 1 } G_j n, \ \ \text{and}\ \ \overline{w}_j = \max_{w\in\mathbb W} G_jw\;.
\]
Regarding the objective function, one can, for example, solve the problem in two 
stages. First, one can maximize the domain, for example via $\rho(d,z) = \sum_{i = 1}^   {{\sf f}_1} z_i$, fix it and then minimize the corresponding performance bound with 
$\rho(d,z) = d$ for the given domain. Alternatively, one can combine these conflicting 
objectives $\rho(d,z) = d - \omega \sum_{i = 1}^{{\sf f}_1} z_i$ with weight $\omega>0$. A more thorough discussion of these objectives can be found in~\citet[Sect.~4.2]{Houska2025a}.

Problem~\eqref{eq::DESIGN_PROBLEM} has $1+{\sf v}(n_{\mathrm{u}}+3)+{\sf f}$ variables, where $\nf$ and $\nv$ are the number of facets and vertices of the template. If $\mathbb{X}$ and $\mathbb{U}$ are given by $n_{\mathbb{X}}$ and $n_{\mathbb{U}}$ constraints, respectively,~\eqref{eq::DESIGN_PROBLEM} has ${\sf v}(3+{\sf f}+n_{\mathbb X}+n_{\mathbb U}) + {\sf e}$ constraints, where ${\sf e}$ is the number of edges in the template. The size of the problem is configurable as the template is user-defined.

The statements below follow from Theorem~\ref{thm::CLF}.
\begin{corollary}
\label{cor::DesignProblem}
Let all technical assumptions of Lemma~\ref{lem::ergodic} and Theorem~\ref{thm::CLF} hold. If~\eqref{eq::DESIGN_PROBLEM} admits a feasible point, $(d^\star, u^\star, y^\star, z^\star, \lambda^\star, \kappa^\star)$, then the following statements hold:
\begin{enumerate}
    \item The pair $(M_{z^\star}, d^\star)$ satisfies~\eqref{eq::CLF}, see Theorem~\ref{thm::CLF}.

    \item The domain of $M_{z^\star}$ is given by
    \[
    \mathrm{dom}(M_{z^\star}) =\mathrm{cvx}\left( V_1 z^\star, \ldots, V_{\sf{v}} z^\star \right)\;.
    \]
    It is a robust control invariant set for the nonlinear system. This follows as $M_{z^\star}$ satisfies~\eqref{eq::CLF}.

    \item Let $\theta_1(x), \ldots, \theta_{{\sf v}}(x) \geq 0$ be interpolation weights with
    \[
    x = \sum_{i=1}^{\sf{v}} \theta_i(x) V_i z^\star \quad \text{and} \quad \sum_{i=1}^{\sf{v}} \theta_i(x) = 1,
    \]
    for all $x \in \mathrm{dom}(M_{z^\star})$. Then the feedback law
    \begin{align}
    \label{eq::mu}
    \mu(x) = \sum_{i=1}^{\sf{v}} \theta_i(x) u_i^\star,
    \end{align}
    is feasible and robust, ensuring the closed-loop system satisfies all state and control constraints for all possible uncertainty realizations.

    \item By Lemma~\ref{lem::ergodic}, $d^\star$ provides an upper bound on the ergodic performance of the control law $\mu$ in~\eqref{eq::mu}. 

    \item If $d^\star \leq 0$, then $M_{z^\star}$ is a robust CLF satisfying
    \begin{equation*}
        M_{z^\star}(x) \geq \min_{u \in \mathbb{U}} \max_{w \in \mathbb{W}} L(x,u) + I_{\mathbb{X}}(x)  + M_{z^\star}(f(x,u) + w) 
    \end{equation*}
    for all $x \in \mathbb{R}^{n_{\mathrm{x}}}$.
\end{enumerate}
\end{corollary}
Corollary~\ref{cor::DesignProblem} establishes~\eqref{eq::DESIGN_PROBLEM} as an effective tool for designing generalized convex CLF functions, as demonstrated by the numerical case study in the next section.

\section{Numerical Case Study}
% \section{Numerical Case Studies}
\label{sec::case_studies}

We consider a robust optimal control problem 
for a Runge-Kutta (order 4) discretization of the Van der Pol oscillator 
\begin{equation*}
% \begin{alignedat}{2}
% &\min_{u}  \ &&\int^{\infty}_{t=0} \frac{1}{2}(x_{2}^2 + u^2) \,{\rm d}x\\
% &\!\text{s.t.} && \left\{
% \begin{aligned}
% &\forall t\in [0,\infty) \\
% &\dot x_1(t) = x_2(t) \\
% &\dot x_2(t) = -x_1(t) - \frac{1}{2} x_2(t) (1 - x_1(t)^2) + x_1(t)u(t)\\
% & x(t)\in\left[-\frac{5}{2},\frac{5}{2}\right]^2, \ u(t) \in \left[-\frac{7}{2},\frac{7}{2}\right]
% \end{aligned}\right.
% \end{alignedat}
\begin{aligned}
&\dot x_1(t) = x_2(t) + w_1(t) , \\
&\dot x_2(t) = -x_1(t) - \frac{1}{2} x_2(t) (1 - x_1(t)^2) + x_1(t)u(t) + w_2(t), \\
% x^+_{1} &= x_{1} + \tau x_{2}\\
% x^+_{2} &= x_{2} - \tau\left(x_1 + \frac{1}{2}\left(1-x_{1}^2\right)-x_1u\right)
\end{aligned}
\end{equation*}
with step size $\tau = 0.05$ for piecewise constant $u$ and $w$ and
\begin{eqnarray*}
\begin{array}{rclrcl}
L(x,u) &=& \frac{\tau}{2}(x_2^2 + u^2), \quad & \mathbb X &=& \left\{ x\in\mathbb R^2 \,\middle| \,  x^\intercal x \leq 9 \right\}, \\[0.1cm]
\mathbb U &=& \left[-2,2\right], \quad & \mathbb W &=& [-0.005,0.005]^2.
\end{array}
\end{eqnarray*}
A configuration triple with $\sf f_1 = 48$, $\sf f_2 =265$, $\sf v = 576$, and $\sf e = 840$ was constructed as a discretization of a 2-hemisphere, adapting the method
of~\citet[Section~V]{Houska2024a} such that $h_2 < 0$ is satisfied. 

Since the right-hand side of the dynamics is twice continuously differentiable, 
Assumption~\ref{ass::f} is satisfied with $\alpha =2$ and $\gamma \geq 0.05$---obtained by bounding the spectral norm of the Hessians on $\mathbb X \times \mathbb U$. Problem~\eqref{eq::DESIGN_PROBLEM} was solved in 
two stages, as explained in the previous section.

\begin{figure}[h!]
    \centering
        \includegraphics[width=0.9\linewidth]{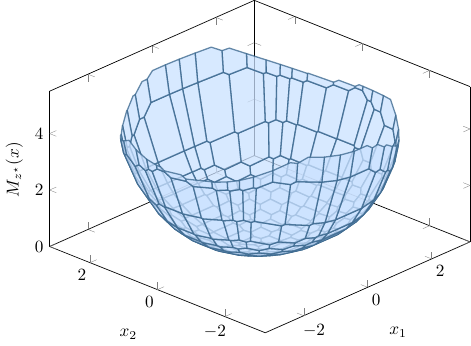}
        \caption{\label{fig::LCLF}A Piecewise Affine Generalized CLF $M_{z^\star}$.}
\end{figure}

\begin{figure}[h!]
    \centering
        \includegraphics[width=0.85\linewidth]{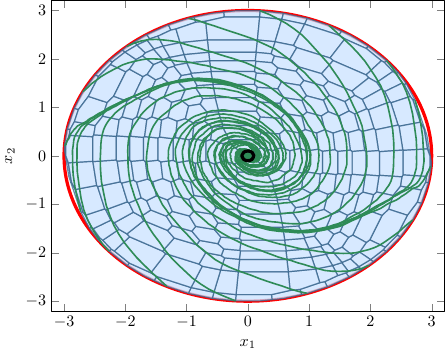}
        \caption{\label{fig::partition}Partition of $\operatorname{dom}(M_{z^\star})$ and closed-loop trajectories.}
\end{figure}

\begin{figure}[h!]
    \centering
        \includegraphics[width=0.9\linewidth]{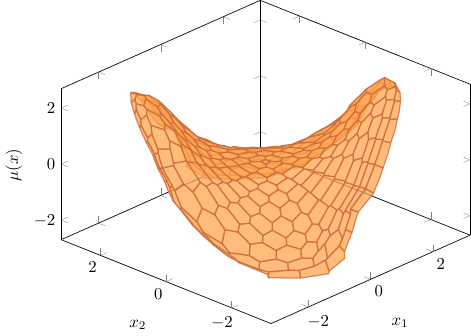}
        \caption{\label{fig::mu}Piecewise affine feedback $\mu$ interpolated from $u^\star$.}
\end{figure}

Figure~\ref{fig::LCLF} shows the generalized CLF $M_{z^\star}$ with 
supply rate $(L-d^\star)$ with $d^\star = 0.1$. This CLF 
induces a partition of the domain 
$\{ x\in\mathbb{R}^2 \,|\, G_{1} x \leq z^\star_{1} \}$ into 265 polytopic regions
shown in Fig.~\ref{fig::partition} together with the state constraint set $\mathbb X$ and selected closed-loop trajectories under random disturbances chosen at each step from the vertices of $\mathbb W$ and the piecewise affine feedback law obtained by least-squares interpolation of the vertex-control inputs $u_1^\star,\ldots,u_{\sf v}^\star$
over the polytopic partition, as shown in Fig.~\ref{fig::mu}. The closed-loop trajectories converge to a small robust control-invariant set near the origin, as illustrated by the black contour in Fig.~\ref{fig::partition}.

\section{Conclusion}
\label{sec::conclusions}

This paper has presented a nonlinear optimization-based framework for synthesizing feedback controllers for constrained, uncertain, nonlinear discrete-time systems. The framework is based on configuration-constrained convex generalized control Lyapunov functions (CLFs) with predetermined storage complexity. They enable the generation of feedback laws with adjustable performance and evaluation times by selecting an appropriate configuration template. Theoretical robustness and performance guarantees have been established by proving that the optimized CLF and its drift parameter satisfy a generalized Lyapunov dissipation inequality, as formalized in Theorem~\ref{thm::CLF} and summarized in Corollary~\ref{cor::DesignProblem}. A numerical case study on a Van der Pol oscillator illustrates the framework's practical applicability. Future research will focus on extending the theoretical results, exploring applications in economic model predictive control, and conducting numerical studies on larger-scale systems.

\small
\bibliography{ifacconf}

\end{document}